\documentclass{amsart}
\usepackage{amssymb,amsmath,amsfonts,latexsym,amscd}
\usepackage{pb-diagram}

\usepackage[margin=1.25in]{geometry}

\usepackage{setspace}
\numberwithin{equation}{section}

\newcommand{\tr}{{\operatorname{tr}}}

\newcommand{\Ind}{{\operatorname{Ind}}}

\newcommand{\Tr}{{\operatorname{Tr}}}

\newcommand{\complex}{{\mathbb C}}

\newcommand{\DD}{{\raise 0.09em \hbox{/}} \kern -.58em {\partial}}
\newcommand{\C}{{\mathbb{C}}}
\newcommand{\R}{{\mathbb{R}}}

\newcommand{\Q}{{\mathbb{Q}}}

\newcommand{\x}{{\times}}

\newcommand{\cala}{{\mathcal A}}

\newcommand{\calc}{{\mathcal C}}

\newcommand{\cale}{{\mathcal E}}
\newcommand{\calf}{{\mathcal F}}
\newcommand{\calg}{{\mathcal G}}
\newcommand{\calh}{{\mathcal H}}

\newcommand{\calp}{{\mathcal P}}

\newcommand{\calr}{{\mathcal R}}

\newcommand{\calz}{{\mathcal Z}}

\theoremstyle{plain}
        \newtheorem{theorem}{Theorem}[section]
        \newtheorem{lemma}[theorem]{Lemma}
        \newtheorem{proposition}[theorem]{Proposition}
        \newtheorem{corollary}[theorem]{Corollary}

\theoremstyle{definition}
        
        \newtheorem{remark}[theorem]{Remark}
        \newtheorem{example}[theorem]{Example}

\title{Hopf algebroids and secondary characteristic classes}

\author{Jerome Kaminker}
\address{Department of Mathematical Sciences\\IUPUI\\Indianapolis, IN}
\curraddr{Department of Mathematics\\UC Davis\\Davis, CA, 95616}
\email{kaminker@math.ucdavis.edu}
\author{Xiang Tang}
\address{Department of Mathematics\\Washington University\\St. Louis, MO, 63130}
\email{xtang@math.wustl.edu}
\thanks{The second author thanks NSF for support}
\date{\today}

\begin{document}
\mbox{ } \vspace{-15mm}\\
\begin{abstract}
  We study a Hopf algebroid, $\calh$, naturally associated to the
  groupoid $U_n^\delta\ltimes U_n$. We show that classes in the Hopf
  cyclic cohomology of $\calh$ can be used to define secondary
  characteristic classes of trivialized flat $U_n$-bundles. For
  example, there is a cyclic class which corresponds to the universal
  transgressed Chern character and which gives rise to the continuous
  part of the $\rho$-invariant of Atiyah-Patodi-Singer. Moreover,
  these cyclic classes are shown to extend to the K-theory of the
  associated $C^{*}$-algebra.  This point of view gives leads to
  homotopy invariance results for certain characteristic numbers.  In
  particular, we define a subgroup of the cohomology of a group
  analogous to the Gelfand-Fuchs classes described by Connes,
  \cite{connes:transverse}, and show that the higher signatures associated to
  them are homotopy invariant.
\end{abstract}
\maketitle \tableofcontents

\section{Introduction}

We will study some new cyclic classes which forms a subgroup of the
cohomology of a discrete group. In common with the Gelfand-Fuchs
classes used by Connes in \cite{connes:transverse}, they extend to
pair with the K-theory of a certain algebra and hence yield homotopy
invariance results for some ``higher signatures''.    The original
goal was to study the Connes-Moscovici index theory for hypoelliptic
transverse signature operators, but in the much easier case when the
foliation is Riemannian.  The spectral triple in that case was
worked out by Wu
 and one needed to find the appropriate adaptation of the Hopf
algebra techniques of \cite{co-mo:gelfand-fuks}.  As a step in that direction
we considered a very special class of Riemannian--foliated flat
$U_{n}$ bundles. We will use a Hopf algebroid to define characteristic
classes and relate them to several earlier constructions.  The
transverse signature operator in this case comes from the signature
operator on the compact Lie group $U_{n}$ and its Chern character is
in the periodic cyclic cohomology of our Hopf algebroid.  However, in
the present paper we will look at different aspects of these classes.

The Hopf algebroid, which we denote $\calh( U_{n}^{\delta}\ltimes
U_n)$, is one associated to the \'etale groupoid, $
U_{n}^{\delta}\ltimes U_n$. We compute its cyclic theory as
introduced by Connes and Moscovici, \cite{co-mo:curved}, and see
that it is isomorphic to the cyclic theory of
$C_{c}^{\infty}(U_{n}^{\delta}\ltimes U_n)$.  Moreover, it has a
convenient description via a double complex.

The double complex has a subcomplex corresponding to invariant forms
which we investigate further.  These classes have two important
properties.  The first is that the associated cocycles extend to
$K_{\bullet}(C(U_{n}) \rtimes U_{n}^{\delta})  $, the K-theory of
the reduced $C^{*}$-algebra.  This yields homotopy invariance
results.  The second is that they can be expressed in terms of
transgressed classes on the base.  A particular case yields the
transgressed Chern character, hence yields the homotopy invariance
of the (continuous part of) the $\rho$-invariant of
Atiyah-Patodi-Singer.

Using the ideas in Connes, \cite{connes:transverse}, we also consider the
subgroup of the cohomology of a discrete group which is determined by
these classes.  All of the classes in this subgroup satisfy the
Novikov conjecture.  It would be interesting to know how much of the
cohomology of a groups is spanned by these classes and Gelfand-Fuchs
classes.  This is analogous to results on low-dimensional cohomology
of groups and the Novikov conjecture in the sense that the results do
not depend on any special properties of the groups.

The authors would like to thank Sasha Gorokhovsky for valuable
discussions, and the first author particularly wants to thank Steve
Hurder and Ron Douglas since many of the ideas developed here came
from his earlier collaboration with them.  We also want to thank
Alain Connes for suggesting these directions several years back.

\section{Cyclic cohomology of  Hopf algebroids}
In this section we will review Connes and Moscovici's definition,
\cite{co-mo:curved}, of cyclic cohomology of a Hopf algebroid.  We will
study a special case related to the classifying space for trivialized flat
unitary bundles.
\subsection{Hopf algebroids}
 In \cite{lu:algebroid}, Lu introduced the notion of a Hopf algebroid as
a generalization of a Hopf algebra. Cyclic cohomology of a Hopf
algebroid was developed by Connes and Moscovici,
\cite{co-mo:curved}, in their study of transverse index theory in
the non-flat case.

Let $A$ and $B$ be unital algebras.
A bialgebroid structure on $A$,
over $B$, consists of the following data.
\begin{enumerate}
\item[i)]  An algebra homomorphism $\alpha:B\to A$ called the
  \emph{source map} and an algebra anti-homomorphism $\beta:B\to A$ called the
\emph{target map}, satisfying $\alpha(a)\beta(b)=\beta(b)\alpha(a)$,
for all $a,b\in A$. \vskip.1in

For the next part of the definition let $A\otimes_B A$ be the
quotient of $A\otimes A$ by the right $A\otimes A$ ideal generated by
$\beta(a)\otimes1-1\otimes \alpha(a)$ for all   $a\in A$.

\item[ii)]   A $B$-$B$
bimodule map $\Delta:
  A\to A\otimes_B A $, called the \emph{coproduct}, satisfying
    \begin{enumerate}
    \item $\Delta(1)=1\otimes 1$;
    \item $(\Delta\otimes_B Id)\Delta=(Id\otimes_B \Delta)\Delta:A\to
      A\otimes_B A\otimes_B A$,
    \item $\Delta(a)(\beta(b)\otimes1-1\otimes\alpha(b))=0$, for $a\in
      A$,  $b\in B$, %(Here one uses the fact that
      %$A\otimes_B A$ is a right $A\otimes A$ module.)
    \item $\Delta(a_1a_2)=\Delta(a_1)\Delta(a_2)$, for $a_1, a_2\in
      A$. % (The product on the right side is
      %well defined by (c).)
    \end{enumerate}
\item  A $B$-$B$ bimodule map $\epsilon:A\to B$, called the
  \emph{counit},  satisfying
    \begin{enumerate}
    \item $\epsilon(1)=1$;
    \item $\ker \epsilon$ is a left $ A$ ideal;
    \item $(\epsilon\otimes_B id)\Delta=(Id\otimes
      \epsilon)\Delta=Id:A\to A$ %(We identify $B\otimes_B
      %A=A=A\otimes_B B$);
    \end{enumerate}
\end{enumerate}

A \emph{Hopf algebroid} is a bialgebroid $A$, over $B$, which admits an
algebra anti-isomorphism $S: A\to A$ satisfying
\begin{enumerate}
\item $S\circ \beta=\alpha$;
\item $m_A(S\otimes id)\Delta=\beta\epsilon S:A\to A$, with
  $m_A:A\otimes A\to A$ the multiplication on $A$;
\item there is a linear map $\gamma:A\otimes_B A\to A\otimes A $ such
  that
    \begin{enumerate}
    \item If $\pi: A\otimes A\to A\otimes_B A$ is the natural
      projection, $\pi\gamma=Id: A\otimes_B A\to A\otimes_B A$;
    \item $m_A(Id\otimes S)\gamma\Delta= \alpha\epsilon: A\to A$.
    \end{enumerate}
\end{enumerate}

We note that in the above definition  one may allow $A$ and $B$ to be
differential graded algebras and require all of the above maps to be
compatible with the differentials and to be of degree 0. Thus one
would have a differential graded Hopf algebroid. The particular example
of a Hopf algebroid which we will study in the present paper is of
this type and is inspired by
\cite{co-mo:curved}, \cite{gro:secondary}.

\begin{example}
\label{ex:main} We consider a smooth \'etale groupoid
$\calg\rightrightarrows G_0$. Define $A$ to be the algebra of
 differential forms\footnote{We remark that we do not require differential forms to be compactly supported.} on $\calg$, and $B$ to be the algebra of
differential forms on $G_0$. Both $A$ and $B$ are differential
graded (commutative) algebras with the de Rham differential.

We define the source and target map $\alpha,\beta: A\to B$ as the
duals of the source and target maps of the groupoid
$\calg\rightrightarrows G_0$.  Similarly, we define $\Delta$ to be
the dual of the product on $\calg \times_{G_0} \calg \to \calg$.
Note that the algebra of differential forms on $\calg \times_{G_0}
\calg$ is equal to $A\otimes_B A$.  The counit map $\epsilon: A\to
B$ is the embedding $G_0\to \calg$ of the unit space.  It is
straightforward to check that $(A,B, \alpha, \beta, \Delta,
\epsilon)$ is a bialgebroid.

The antipode map $S:A\to A$ is defined to be the dual of the inversion map
 $\iota :\calg\to \calg$.  It is easy to check that $S$ satisfies properties (1) and (2) for
an antipode of a Hopf algebroid. In the case that
$\calg\rightrightarrows G_0$ is a transformation groupoid from a
discrete group $\Gamma$ action on $G_0$, we define $\gamma$ to be
the dual of the map $\calg\times \calg \to \calg\times _{G_0}\calg$
with $\calg\times \calg\ni (\alpha_1,m_1,\alpha_2,m_2)\mapsto
(\alpha_1, m_1, \alpha_2, \alpha_2^{-1}(m_2))\in
\calg\times_{G_0}\calg$. It is not difficult to check that $(A, B,
\alpha, \beta, \Delta, \epsilon, S, \gamma)$ is a Hopf algebroid.
Since   $A$ and $B$ are topological algebras, we will want to work
in that category so we will use the projective tensor product.  This
Hopf algebroid will be referred to as $\calh ( \calg)$.  In general,
we do not know a natural construction of the linear map $\gamma$
required in the above property (3) of an antipode. We refer to
\cite{mu:hopf} for a further discussion about Hopf algebroid
associated to a smooth \'etale groupoid.
\end{example}
\subsection{Cyclic cohomology of a Hopf algebroid and Hopf algebroid  modules}
\label{sec:def-cohomology} For completeness, we will review Connes
and Moscovici's definition of cyclic cohomology of a Hopf
algebroid.  The following construction works for a bialgebroid
with a \emph{ twisted}  antipode $S:A\to A$.  That is, one has an algebra
anti-isomorphism such that
\[
S^2=id,\ \ \ \ \ S\beta=\alpha,\ \ \ \ \ m_A(S\otimes_B
Id)\Delta=\beta\epsilon S:A\to A,
\]
and
\begin{equation}\label{eq:antipode}
S(a^{(1)})^{(1)}a^{(2)}\otimes_B S(a^{(1)})^{(2)}=1\otimes_B S(a).
\end{equation}
In the above formula we have used Sweedler's notation for the
coproduct $\Delta(a)=a^{(1)}\otimes_B a^{(2)}$. We remark that a
bialgebroid with this twisted version of an antipode is called a
para-Hopf algebroid in \cite{kh-ra:hopf-cyclic}.
\begin{proposition}\label{prop:groupoid}
The bialgebroid $\calh(\calg)=(A,B,\alpha,\beta, \Delta, \epsilon)$ with the antipode $S$ in Example \ref{ex:main} for a general \'etale groupoids forms a para-Hopf algebroid.
\end{proposition}
\begin{proof}It is sufficient to check Equation (\ref{eq:antipode}). Let $a$ be an element in $A$. For $x,y\in \calg$ two composable arrows, we notice that
\[
\Delta(a)(x,y)=a(xy),\ \ \ (S\otimes Id)\Delta(a)(x,y)=a(x^{-1}y).
\]
Therefore, one computes
\[
\begin{split}
(\Delta S\otimes Id)\Delta(a)(x,y,z)&=a(y^{-1}x^{-1}z),\\
S(a^{(1)})^{(1)}a^{(2)}\otimes_B
S(a^{(1)})^{(2)}(x,y)&=a(y^{-1}x^{-1}x)=a(y^{-1})=1\otimes_B S(a).
\end{split}
\]
\end{proof}
Let $\Lambda$ be the cyclic category. We construct a cyclic module
$A^\natural$ for a para-Hopf algebroid $(A, B, \alpha, \beta,
\Delta,\epsilon, S)$.

Define
\[
C^0=B,\ \ \ \ \ C^n=\underbrace{A\otimes _B A\otimes _B\cdots \otimes_B A}_n,\ n\geq 1.
\]

Faces and degeneracy operators are defined as follows:
\[
\begin{split}
\delta_0(a^1\otimes_B\cdots \otimes_B a^{n-1})&=1\otimes_B a^1\otimes_B \cdots \otimes_B a^{n-1}\\
\delta_i(a^1\otimes_B\cdots \otimes_B a^{n-1})&=a^1\otimes_B\cdots \otimes_B \Delta a^i\otimes_B
\cdots \otimes_B a^{n-1},\ \ \ \ \ 1\leq i\leq n-1\\
\delta_n(a^1\otimes_B \cdots \otimes_B a^{n-1})&=a^1\otimes_B\cdots \otimes_Ba^{n-1}\otimes_B1;\\
\sigma_i(a^1\otimes_B\cdots \otimes_B a^{n+1})&=a^1\otimes_B\cdots
\otimes_B a^i\otimes_B \epsilon(a^{i+1})\otimes_Ba^{i+2}\otimes_B
\cdots \otimes_B a^{n+1}.
\end{split}
\]

The cyclic operators are given by
\[
\tau_n(a^1\otimes_B\cdots \otimes_B
a^n)=(\Delta^{n-1}S(a^1))(a^2\otimes \cdots a^n\otimes 1).
\].
% where on the right side of the above equation, we have used the
% property that the right multiplication between  $A\otimes \cdots
% \otimes A$ and  $A\otimes_B\cdots \otimes_B A$ is well defined.

We define the Hopf cyclic cohomology of $A$ to be the cyclic
cohomology of $A^\natural$. We remark that this definition extends
naturally to differential graded Hopf algebroids since all the maps
above are maps between differential graded vector spaces and the
differential is of degree 1. In the following, we briefly review
Gorokhovsky's construction \cite{gro:secondary} of characteristic
class map for a differential graded Hopf algebroid action. We remark
that Gorokhovsky only studied characteristic map for a differential
graded algebra action. The same construction works for a general
para-Hopf algebroid action. In the following we will only state the
results and their proofs are actually identical to those in
\cite{gro:secondary}.

A differential graded algebra $M$ is equipped with a Hopf algebroid
$(A, B, \alpha, \beta, \Delta, \epsilon, S)$ action, if there is a
differential graded algebra morphism $\rho:B\to M$ and an action
$\lambda:A\otimes_B M\to M$ satisfying (1) $\lambda(a_1a_2,
m)=\lambda(a_1, \lambda(a_2, m)$ for $a_1, a_2\in A, m\in M$; (2)
$deg(\lambda(a,m))=deg(a)+deg(m)$ and $d(\lambda(a,m))=\lambda(da,
m)+(-1)^{deg(a)}\lambda(a,dm)$ for $a\in A, m\in M$; (3)$\lambda(1,
m)=m$ for $m\in M$; (4)
$\lambda(a,m_1m_2)=\sum_{i}(-1)^{deg(m_1)deg(a_2^i)}\lambda(a_1^i,
m_1)\lambda(a_2^i, m_2)$, for $a\in A$ and $\Delta(a)=\sum_i
a_1^i\otimes a_2^i\in A\otimes _B A$; and when $M$ has a unit, we
require (5) $\lambda(a,1)=\lambda(\epsilon(a),1)=\rho(\epsilon(a))$
for $a\in A$.

We assume that there is a trace $\tau$ on $M$ such that
\[
\tau(a(m)n)=\tau(mS(a)(n)), \ \ \ \ \ \ \ \forall m, n \in M, \ \forall a\in A.
\]

We define a cochain map $c$ from cyclic cochain complex of $A$ to
cyclic cochain complex of the differential graded algebra $M$, i.e.
for $\Phi=\sum_i a^i_1\otimes \cdots \otimes a^i_q$ an element in $
C^{q}$,
\[
c(\Phi)(m_0, \cdots, m_q)=\tau(\sum_i m_0\lambda(a^i_1,m_1)\cdots
\lambda(a^i_q,m_q)),\ \ \ \ \ \ m_i\in M.
\]
We call a cyclic cochain on $M$ in the image of the form $c(\Phi)$ a
differentiable cochain, and denote the space of differentiable
$p$-cochains by $C^p_d(M)$. Then we can define differentiable cyclic
cohomology of $M$, $HC^*_d(M)$, to be the cyclic cohomology of
$C^\natural_d(M)$. One can also define the periodic differentiable
cyclic cohomology, $HP^*_d(M)$, by similar means. One has  the
following natural map
\begin{eqnarray*}
c:HC^*(A)\rightarrow HP^*_d(M)\rightarrow HP^*(M).
\end{eqnarray*}

For later applications, we assume that the trace $\tau$ is of weight
$q$, i.e. $\tau(m)=0$, if $deg(m)\ne q$. Then we easily see that a
cochain in $A^\natural$ with total degree greater than $q$ is mapped
to zero by $c$. Define $F^lA^\natural=\{a_1\otimes \cdots a_j|
deg(a_1)+\cdots deg(a_j)\geq l\}$, and
$A^\natural_q=A^\natural/F^{q+1}A^\natural$. Then $A^\natural_q$ is
again a cyclic module with cyclic cohomology $HC^*(A)_q$ and $
HP^*(A)_q$. And the characteristic map $c$ induces a map $c:
HC^*(A)_q\to HP^*(M)$.

Let $M^0$ be the degree 0 part of $M$. We can map cyclic cochains on
$M$ to $M^0$. And we have the following characteristic map analogous
to \cite{gro:secondary}[Theorem 8],
\[
c:HC^{*}(A)_q\to HP^{*-q}(M^0).
\]

For an \'etale groupoid $\calg$, we consider the space $M$ of
compactly supported differential forms on $\calg$ with the
convolution product. As $\calg$ is an \'etale groupoid, the unit
space $G_0$ is embedded in $\calg$ as an open submanifold. This
embedding map between spaces defines an algebra homomorphism $\rho$
from $B=\Omega^*(G_0)$ to $M$. For $a\in A=\Omega^*(\calg)$ and
$m\in M=\Omega^*_c(\calg)$, define $\lambda(a,m)=am\in
\Omega^*_c(\calg)$ by viewing both $a$ and $m$ as elements of
$\calh(\calg)$ and using the product on $\calh(\calg)$. It is a
direct check that $\rho$ and $ \lambda$ defines an action of
$\calh(\calg)$ on $M$. Furthermore, we assume that there is a
$\calg$ invariant volume form $\Omega$ on $G_0$. The integration on
$G_0$ with respect to $\Omega$ defines a trace $\tau_\Omega$ on $M$
of weight $\dim(G_0)$ compatible with the antipode $S$. Hence as
$M^0=C^\infty_c(\calg)$, we have a characteristic map $c_\Omega:
HC^*(\calh(\calg))_{\dim(G_0)}\to
HP^{*-\dim(G_0)}(C^\infty_c(\calg))$.

\subsection{Hopf cyclic cohomology of the Hopf algebroid $\calh(\calg)$}
In this section we will compute the Hopf cyclic cohomology of the
para-Hopf algebroid $\calh(\calg)$ introduced in Proposition
\ref{prop:groupoid} for a smooth \'etale groupoid.

To formulate the results, we will use the following constructions.
Let $B\calg^n= \underbrace{\calg\times_{G_0} \cdots\times_{G_0}
\calg}_n$, where $\calg \x_{G_0} \calg$ is the fiber product with
respect to the maps $t:\calg \to G_0$ and $s:\calg\to G_0$.  This
can be given the structure of a simplicial manifold.  Let
$\Omega^{*}(B\calg^{*}) $ be the double complex of sheaves of
differential forms on the simplicial manifold $B\calg^*$,
\cite{du:simplicial}. $\Omega^*(B\calg^*)$ has a standard cyclic
structure, (described below). As $(\Omega^*,d)$ is a resolution of
$\C$, we denote the corresponding cyclic cohomology of
$\Omega^{*}(B\calg^{*}) $ by $HC^*(\calg, \C)$ and $HP^*(\calg,\C)$.
One can now state the following result.
\begin{theorem}
\label{thm:cohomology-hopf}
\[
HC^*(\calh(\calg))_{\dim(G_0)}=HC^{*}(\calg; \complex)\ \
\text{and}\ \ HP^*(\calh(\calg))_{\dim(G_0)}=HP^{*}(\calg;
\complex).
% HP^\bullet(\calh(\calg))&=HP^\bullet(\calg; \complex )=\bigoplus_{\bullet+2k,
% k\in \mathbb{Z}}H^\bullet(\calg; \complex),
\]
% where $H^\bullet(\calg;\complex)$ is the groupoid
% sheaf cohomology of the sheaf $\complex$.
\end{theorem}
\noindent{\bf Proof:}
As a first step one notes that an
% We consider the space $B\calg_n=
% \underbrace{\calg\times_{M} \cdots\times_M \calg}_n$.
% We identify
 $m$-cochain on $\calh(\calg)$ of total degree less than or equal to $\dim(G_0)$
can be identified with a smooth section of the sheaf of differential
forms on $B\calg^m$, $m\geq 0$. This identification respects the
cyclic simplicial structures on $\Omega^{*}(B\calg^{*}) $
% $\Gamma^\infty(B\calg_n; \Omega^\bullet)$
and $\calh(\calg)^{\natural}$. Indeed, following the notation in
\cite{cr:cyclic}, where $(a\ |\  g_1, \cdots, g_n)$ with $a\in
\Omega^*(G_0)|_{t(g_1)}$, $g_i\in \calg$ and $s(g_i)=t(g_{i+1}),
i=1,\cdots, n-1$ denotes elements of $\Omega^{*}(B\calg^{*}) $
% $\Gamma^\infty(B\calg;
% \Omega^*)$
, the cyclic simplicial structure on $\Omega^{*}(B\calg^{*})$ is
given by
\[
\delta_i(a\ |\  g_1, \cdots, g_n)=\left\{\begin{array}{ll}(ag_1\ |\
g_2, \cdots, g_n)& i=0\\ (a\ |\  g_1, \cdots, g_ig_{i+1}, \cdots,
g_n)& 1\leq i\leq k-1\\ (a\ |\  g_1, \cdots, g_{n-1})&
i=n\end{array}\right.,
\]
and
\[
\sigma_i(a\ |\  g_1, \cdots, g_n)=(a\ |\  \cdots, g_{i-1}, 1, g_{i},
\cdots, g_n),
\]
and
\[
t(a\ |\  g_1, \cdots, g_n)=(ag_1\cdots g_n\ |\  (g_1g_2\cdots
g_n)^{-1}, g_1, \cdots, g_{n-1}).
\]
In the above,  $ag$  means the right translation of $a$ by $g\in \calg$.

It is easy to check from this identification that the Hopf cyclic
cohomology of $\calh(\calg)$ is equal to the cyclic cohomology of
$\Omega^{*}(B\calg^{*}) $.
% $\Gamma^\infty(B\calg; \Omega^*)$
$\Box$

Just as the cohomology of a discrete group is isomorphic to the
cohomology of its classifying space, there is an analogous result
for \'etale groupoids, \cite{moerdijk}.  It states that the groupoid
sheaf cohomology of the sheaf $\complex$ on $\calg$, denoted
$H^*(\calg; \complex)$, is isomorphic to the cohomology of its
classifying space $B\calg$.

Since $\Omega^*$ is a projective resolution of the sheaf $\complex$
on $\calg$,  we obtain the following corollary.
\begin{corollary}
\[
\begin{split}
HC^*(\calh(\calg))_{\dim(G_0)}&=\bigoplus_{k\in
\mathbb{Z}_+}H^{*-2k}(B\calg; \complex),\\
HP^*(\calh(\calg))_{\dim(G_0)}&=\bigoplus_{k\in
\mathbb{Z}}H^{*+2k}(B\calg; \complex).
\end{split}
\]
\end{corollary}

We apply Theorem \ref{thm:cohomology-hopf} in the following special case.

Let $U_n$ be the compact Lie group of $n\times n$ unitary matrices
and let $U_n^\delta$ be the same group, but equipped with the
discrete topology.  Left multiplication of $U_n^\delta$ on $U_n$
defines a transformation groupoid $\calg=U_n^\delta\ltimes
U_n\rightrightarrows U_n$. We will compute the cyclic cohomology of
its groupoid algebra, $C^\infty_c(U_n^\delta\ltimes U_n)$. The
following result was first proved in \cite{br:spectral}, and we
derive it as an application from \cite{cr:cyclic}.

\begin{proposition}
\label{prop:cohomology-groupoid}
\[
HP^*(C_c^\infty(U_n^\delta\ltimes U_n))=\bigoplus_{ k\in
\mathbb{Z}}H^{*+2k+\dim(U_n)}(U_n^\delta\ltimes U_n, \C).
\]
\end{proposition}
\noindent{\bf Proof:} Because $U_n^\delta$ acts on $U_n$ freely, the
cyclic groupoid $\calz$ introduced in \cite{cr:cyclic} associated to
the \'etale groupoid $U_n^\delta\ltimes U_n$ is equal to itself. Let
$or$ be the orientation sheaf on $U_n^\delta\ltimes U_n$. Then
\cite{cr:cyclic}[4.14] proves $HP^*(C^\infty_c(U_n^\delta\ltimes
U_n))=\bigoplus_{k\in
\mathbb{Z}}H^{*+2k+\dim(U_n)}(U_n^\delta\ltimes U_n, or)$. Using the
left invariant volume form $\Omega$ on $U_n$, one can identify the
orientation sheaf $or$ and $\C$, and the statement of this
proposition follows. $\Box$

Note that, in the present situation, the groupoid sheaf cohomology,
$H^*(U_n^\delta\ltimes U_n, \complex)$, is by definition the
cohomology of the total complex of
\[
\Gamma( \big(\underbrace{U_n^\delta\times\cdots \times
U_n^\delta}_{p}\big)\times_{U_n^\delta}U_n, \Omega^q(U_n)),
\]
and the result is equal to the simplicial de Rham cohomology
$H^*(EU_n^\delta\times_{U_n^\delta}U_n)$.  We record this fact for
later use.
\begin{corollary}
  \label{cohomology}
The groupoid sheaf cohomology, $H^*(U_n^\delta\ltimes U_n,
\complex)$ is isomorphic to the de Rham cohomology of the simplicial
manifold $H^*(EU_n^\delta\times _{U_n^\delta} U_n)$.

\end{corollary}

Let $\overline{BU}_{n} $ denote the geometric realization of the
simplicial manifold $EU_n^\delta\times _{U_n^\delta} U_n$.  Recall
that there is a Connes' map $\phi: \oplus_{k}
H^{*+2k}(\overline{BU}_{n}) \to
HP^{*-\dim(U_n)}(C_{c}^{\infty}( U_{n}^{\delta}\ltimes U_n))$.
This map actually agrees with $c:\oplus_k
H^{*+2k}(\overline{BU}_n)=HP^*(\calh(U_n^\delta\ltimes
U_n))_{\dim(U_n)}\to HP^{*-\dim(U_n)}(C^\infty_c(U_n^\delta\ltimes
U_n))$ by the same argument as \cite{gro:secondary}. We observe that
the composition of the map $c$ with the identification
$e:HP^*(C^\infty_c( U_n^\delta\ltimes
U_n^\delta))=HP^{*+\dim(U_n)}(U_n^\delta\ltimes U_n, \C)$ is the
identity map on $H^*(U_n^\delta\ltimes U_n, \C)$ as $c$ maps a
differential form on $U_n^\delta\ltimes U_n$ to a differential
current on $U_n^\delta\ltimes U_n$ to pair with cyclic chains on
$C_c^\infty(U_n^\delta \ltimes U_n)$ and $e$ maps vice versa. In
summary, we have obtained the following corollary.
\begin{corollary}
\label{phi} The map $$\phi: \oplus_k
H^{*+2k}(\overline{BU}_{n})=HP^*(\calh(U_n^\delta\ltimes
U_n))_{\dim(U_n)}\to HP^{*-\dim(U_n)}(C_{c}^{\infty}(
U_{n}^{\delta}\ltimes U_n)),$$ is an isomorphism.
\end{corollary}

C. Lazarov and J. Pasternack have introduced secondary
characteristic classes for Riemannian foliations with trivial normal
bundle. We will  relate our constructions to theirs,
\cite{la-pa:secondary}.  Let $\overline {BR\Gamma}_{q}$ denote the
classifying space for Riemannian Haefliger structures of codimension
$q$ with trivial normal bundle. We remark that the dimension of
$U_n$ is $n^2$.

\begin{proposition}
\label{prop:pasternack} There is a canonical Riemannian Haefliger
structure on $EU_n^\delta\times _{U_n^\delta} U_n$ with trivial
normal bundle which induces a map $ EU_n^\delta\times _{U_n^\delta}
U_n\to \overline {BR\Gamma}_{n^{2}}$.

% There exists a natural map
% \[
% \rho: H^*(RW_{\frac{n(3n+1)}{2}})\to
% H^*(EU_n^\delta\times _{U_n^\delta} U_n).
% \]

\end{proposition}

%\begin{aside}
  % \bf{I think the proof below essentially does this and I will fill in
%     the construction of the necessary Haefliger cocycle.  It is an
%     interesting question whether we are getting some new nontrivial
%     classes in $H^{*}(\overline {BR\Gamma}_{n^{2}})$.  By this I mean
%     that an index theorem might show that the image of some class in
%     $H^{*}(\overline {BR\Gamma}_{n^{2}})$ } gives a cyclic class which
%   pairs non trivially.
%\end{aside}
\noindent{\bf Proof:}
For the existence of a Riemannian Haefliger structure cover
$EU_{n}^{\delta}$ by small open sets $V_{\alpha}$ which are disjoint
from all their translates.  Choose an invariant Riemannian metric on
$U_{n}$ and let $\{W_{\beta}\}$ be a cover by charts.  Then $pr_{1}:
V_{\alpha} \x W_{\beta} \to \R^{n^{2}}$ descends to the quotient by
$U_{n}^{\delta}$ and provides the Haefliger structure.  We must
show that it's normal bundle is trivial.
% We notice that $U_n\times
% _{U_n^\delta}EU_n^\delta$ is foliated by leaves of the form
% $pt\times EU_n^\delta$. We denote this foliation by $\calf_{U_n}$.
% In the following we prove that the normal bundle of this foliation
% is actually trivial.
% Consider the projection $pr_{1}:U_n\times
% EU_n^\delta\to U_n$  onto the first component.
Note first that the tangent bundle $TU_n$ is a trivialized by using
left translation.  Then $pr_{1}^*(TU_n)$  is again trivial. Because
$U_n^\delta$ acts trivially on $TU_n$, the quotient
$\widetilde{TU_n}$ of the bundle $p_{1}^*TU_n$ on $EU_n^\delta\times
_{U_n^\delta} U_n$ is again trivial. But the coordinate
transformations for
 $\widetilde{TU_n}$ are obtained from the differentials of those
 associated to the Haefliger structure, so the conclusion holds. $\Box$

% We have shown that on $EU_n^\delta\times _{U_n^\delta} U_n$ there is
% foliation $\calf_{U_n}$ which has a trivial normal bundle
% $\tilde{T}_{U_n}$.
By \cite{la-pa:secondary}, there is a
characteristic class map obtained as the following composition.
\begin{equation}
\label{eq:rho} \kappa:H^*(RW_{n^{2}})\to H^{*}(\overline
{BR\Gamma}_{n^{2}}) \to H^*(EU_n^\delta\times _{U_n^\delta} U_n).
\end{equation}

\section{Secondary characteristic classes of trivialized flat bundles}
Let $M$ be a compact smooth manifold and let $\Gamma=\pi_1(M)$ be
the fundamental group of $M$. Let $\widetilde{M}$ be the universal
cover of $M$. Suppose that we are given a finite dimensional unitary
representation, $\alpha: \Gamma\to U_n$. Consider the flat principal
$U_n$ bundle $V=\widetilde{M}\times_\Gamma U_n \to M$. If we assume
that this bundle is trivial with a given trivialization
$\theta:\widetilde{M}\times_\Gamma U_n\to M\times U_n$ then we can
relate it to the Lazarov-Pasternack map.

In this section we will construct a map
$\chi:HC^{*}(\calh(U_n^\delta\ltimes {U_n}))_{n^2} \to H^{*}(V)$.
It will be a composition of several maps.

First, the homomorphism of algebras, $C^\infty_c(\Gamma\ltimes
U_n)\to C^\infty_c(U_n^\delta\ltimes U_n) $ induces a homomorphism
on cyclic cohomology
\begin{equation}
\label{eq:restriction} \nu: HP^*(C^{\infty}_{c}(U_n^\delta\ltimes
U_n))\to HP^*(C_{c}^\infty(\Gamma \ltimes U_n)).
\end{equation}

Next, note that the manifold $V$ is foliated by leaves which are the
images of $\widetilde{M}\times \{g\} $. We denote this foliation by
$\calf$, and the corresponding holonomy groupoid by $\calg$. We
assume that the representation $\alpha$ is faithful, which will
assure that the holonomy groupoid $\calg$ is Hausdorff. Furthermore,
$\calg$ is Morita equivalent to the transformation groupoid $\Gamma
\ltimes U_n$. Therefore the groupoid algebra $C^\infty_c(\calg)$ is
Morita
equivalent to $C^\infty_c(\Gamma \ltimes U_n)$. % More precisely,
% $C_c^\infty(U_n)\rtimes \Gamma$ is embedded into
% $C_c^\infty(\calg)$ as a corner.
Thus, there is an isomorphism on cyclic cohomology, which provides the
second map.
\begin{equation}
\label{eq:iso-morita}\iota: HP^*(C^\infty_c(\Gamma \ltimes U_n)) \to
HP^*(C^\infty_c(\calg)).
\end{equation}

Recall that in \cite{connes:book}[Chapter III, 7.$\gamma$], Connes
constructs a map $\lambda: HP^*(C^\infty_c(\calg))\to H^*(B\calg)$
which is a left inverse to the map $\Phi$ which is important in
higher index theorems. Finally, the foliated structure on $V$
defines a map $V \to B\calg$. The induced map on cohomology will be
denoted $\eta: H^*(B\calg)\to H^*(V)$.

We summarize the above constructions in the following sequence,
\[
\begin{array}{l}
HC^*(\calh(U_n^\delta\ltimes
U_n))_{n^2}\stackrel{c}{\longrightarrow}HP^{*-n^2}(C^\infty_c(U_n^\delta\ltimes
U_n))\stackrel{\nu}{\longrightarrow}
HP^{*-n^2}(C^\infty_c(\Gamma\ltimes U_n))\\
\stackrel{\iota}{\longrightarrow}
HP^{*-n^2}(C^\infty_c(\calg))\stackrel{\lambda}{\longrightarrow}H^{*-n^2}
{(B\calg)}\stackrel{\eta}{\longrightarrow}\oplus_kH^{*+2k-n^2}(V).
\end{array}
\]

We define $$\chi: HC^*(\calh(U_n^\delta\ltimes
U_n))_{n^2}\longrightarrow \oplus_k H^{*+2k-n^2}(V)$$ to be the
above composition and
$$\hat \chi : HC^*(\calh(U_n^\delta\ltimes
U_n))_{n^2}{\longrightarrow} HP^{*-n^2}(C^\infty_c(\Gamma \ltimes
U_n))$$ to be the composition of the first two arrows. We will refer
to $\chi$ as the characteristic map and cohomology classes in the
image of $\chi$ will be viewed as secondary characteristic classes.

These classes are compatible with those introduced by Lazarov and
Pasternack as the next proposition shows.

\begin{proposition}
  \label{prop:hopf-pasternack} Let $\kappa$ and $\tilde \kappa$ be the
  maps obtained using Lazarov-Pasternack \cite{la-pa:secondary}
  (\ref{eq:rho}). Then the following diagram commutes.
\[
\begin{diagram}
\node{H^*(RW_{n^{2}})}
\arrow{e,t,T}{\kappa}\arrow{se,t,T}{\tilde \kappa}\node{H^*(EU_n^\delta\times _{U_n^\delta} U_n)}\arrow{s,l,T}{\chi}\\
\node{}\node{H^*(V)}
\end{diagram}
\]
\end{proposition}

\noindent{\bf Proof:} It is straightforward check that the diagram
is commutative. $\Box$

In the remainder of this section we provide a more detailed
description of the map $\iota: HP^*(C^\infty_c(U_n\rtimes
\Gamma))\to HP^*(C^\infty_c(\calg))$ following ideas from
\cite{connes:book}[III. 4 $\alpha$].
% \begin{aside}
%   \bf{This is where I am.  I think what you have below is correct, but
%   some statements may be off a little.  Consider the Kronecker
%   foliation of $T^{2}$.  Each leaf is an immersed copy of $\R$, but
%   is not equal to $\tilde M /\Gamma = \R / \Z$.  I know you are aware
%   of this, but I want to go through it more carefully and I have no
%   more time tonight.}
% \end{aside}

Consider the holonomy groupoid $\calg\rightrightarrows V$. The space $\calg$
consists of holonomy classes of paths connecting two points on a
leaf.
% Each single leaf of the foliation on $V$ is diffeomorphic to
% $M=\widetilde{M}/\Gamma$.
The set of homotopy classes of paths on a single leaf is equal to
the fundamental groupoid of the leaf, $\widetilde{M}\times
\widetilde{M}/\Gamma$, where $\Gamma$ acts on $\widetilde{M}\times
\widetilde{M}$ diagonally. Since the leaves are all simply
connected, it follows from this that the holonomy groupoid $\calg$
can be identified
 with $ (\widetilde{M}\times
\widetilde{M})\times _\Gamma U_n $.  The structure maps $s,t:
\calg=(\widetilde{M}\times \widetilde{M})\times _\Gamma U_n
\rightrightarrows V= \widetilde{M}\times _{\Gamma} U_n$ and groupoid
operation are defined as follows,
\[
\begin{array}{ll}
s<\tilde{x},\tilde{y},g>\  =\  <\tilde{x},g>,\ \ \ t<\tilde{x},
\tilde{y},g>\
=\  <\tilde{y},g>,\ \ \ \ &\tilde{x}, \tilde{y}\in \widetilde{M},\ g\in U_n\\
<\tilde{x}, \tilde{y},g> \circ <\tilde{y}, \tilde{z},g> \ = \
<\tilde{x}, \tilde{z},g>,\ \ \ \ \  &\tilde{x}, \tilde{y},
\tilde{z}\in \widetilde{M}.
\end{array}
\]
Using this description of $\calg$, we can make the following identifications.

Let $\calr$ be the algebra of infinite matrices $(a_{ij})_{i,n\in
\mathbb{N}}$ with rapid decay property, i.e.
\[
sup_{i,j\in \mathbb{N}}i^kj^l|a_{ij}|<\infty, \ \forall k,l\in \mathbb{N}.
\]
One defines the map
\[
T: C^\infty_c(\calg)\longrightarrow M_N\big(C^\infty_c(\Gamma
\ltimes U_n)\otimes \calr\big).
\]
where $M_N\big(C^\infty_c(\Gamma \ltimes U_n)\otimes \calr\big)$ is
the algebra of $N\times N$ matrices with entries in
$C^\infty_c(\Gamma \ltimes U_n)\otimes \calr$.

Since $M$ is a compact manifold, we can choose a finite open cover
$(U_i)_{i=1,\cdots, N}$, with $\beta_i:U_i\to \widetilde{M}$ a local
smooth section of the projection $\widetilde{M}\to M$, and
$(\varphi_i)_{i=1,\cdots, N}$ a smooth partition of unity
subordinate to the covering $(U_i)_{i=1,\cdots, N}$, with
$\varphi_i^{\frac{1}{2}}$ also smooth functions.

Let $\calr_M$ denote the algebra of smoothing operators on $M$.  For
any $f\in C_c^\infty(\calg)$, define $T(f)\in
M_N\big(C^\infty_c(\Gamma \ltimes U_n)\otimes \calr_M\big)$, where
$\calr_M$ is the algebra of smoothing operators on $C^\infty(M)$ by
\[
T(f)_{i,j}(x,y, \gamma,
g)=\varphi_i^{\frac{1}{2}}(x)\varphi_j^{\frac{1}{2}}(y)f(\beta_i(x),\gamma^{-1}\beta_j(y)
,g),\ \ x,y\in M,\  g\in U_n,\ \gamma\in \Gamma.
\]

It is straightforward to check that, since $f$ is smooth and
compactly supported on $\calg$, $T(g)_{ij}(\gamma,g)$ defines a
smoothing operator on $M$ with finite support on $\Gamma$.  Similar
to \cite{connes:book}[III. 4. $\beta$], $T$ is an algebra
homomorphism and induces a Morita equivalence between
$C^\infty_c(\calg)$ and $C^\infty_c(\Gamma \ltimes U_n)$.

We choose an isomorphism $\varrho$ between $\calr_M$ and $\calr$,
which always exists and is unique up to inner automorphisms. The
composition $I=\varrho\circ T$ defines an algebra homomorphism from
$C^\infty_c(\calg)$ to $M_n(C_c^\infty(\Gamma \ltimes U_n)\otimes
\calr).$

We can use $I$ to define a map $\iota: HP^*(C^\infty_c(\Gamma\ltimes
U_n))\to HP^*(C^\infty_c(\calg))$ as follows.  For $\Psi\in
HP^j(C_c^\infty(\Gamma \ltimes U_n))$,
\[
\iota(\Psi)(a_1, \cdots, a_j)=\Psi\big(I(a_1)_{i_1i_2}^{m_1m_2},
I(a_2)_{i_2i_3}^{m_2m_3},\cdots, I(a_n)_{i_ji_1}^{m_jm_1}\big),
\]
where $i_1, \cdots, i_j$ run from 1 to $\infty$, and $m_1,
\cdots, m_j$ runs from  1 to $N$. This construction can be viewed
as the sharp product on cyclic cohomology between $\Psi$ with the
standard traces on $M_n$ and $\calr$.

\section{Extendability of cocycles}
In this section, we will describe a subgroup of  the cohomology of
$\calh(U_n^\delta\ltimes U_n)$, whose image under $\chi$ consists of
classes extending to pair with the K-theory of the reduced
$C^*$-algebra of $\Gamma \ltimes U_n$. The
 proof is a direct adaptation of Connes' argument in
\cite{connes:transverse}.

To describe the subset, consider the double complex
$\calc^{p,q}:=(\Omega^p(U_n)\times_{U^\delta_n} (U_n^{\delta})^{\times
  (q+1)}, (d, b))$ whose total cohomology is
$H^*(U_n^\delta\ltimes U_n)$. By \cite{br:spectral}, the spectral
sequence of this double complex degenerates at $E_2=H^q(BU^\delta_n,
H^p(U_n))$. We consider a subcomplex $(\Omega^{*}(U_n)^{U_n}, d)=
\Omega^*(U_n)\times_{U^\delta_n}
U^\delta_n=\calc^{*,0}\hookrightarrow \calc^{p,q}$. It is easy to
check that the inclusion is a cochain map and therefore defines a
homomorphism $H^*(U_n)^{U_n}\hookrightarrow H^0(BU_n^\delta,
H^*(U_n)) \subseteq HC^*(\calh( U_n^\delta\ltimes U_n))_{n^2}$.

\begin{theorem}\label{thm:extend}
The linear functionals on $K_*(C_{c}^{\infty}(\Gamma\ltimes U_n)$
defined by elements in $\hat \chi(H^*(U_n)^{U_n}) \subset
 HP^{*-n^2}(C^\infty_c(U_n^\delta\ltimes U_n) )$ extend to linear
functionals on $K_*(C(U_n)\rtimes \Gamma)$, where $C(U_n)\rtimes
\Gamma$ is the reduced crossed product $C^*$-algebra.
\end{theorem}

\begin{proof}
We consider the following sequence of maps
\[
H^*(U_n)^{U_n}\rightarrow HC^*(\calh( U_n^\delta\ltimes
U_n))_{n^2}\rightarrow HP^{*-n^2}(C^\infty_c(U_n^\delta\ltimes
U_n))\rightarrow HP^{*-n^2}(C^\infty_c(\Gamma \ltimes U_n)).
\]

Let $\omega\in \Omega^*(U_n))^{U_n}$ be a representative of a
cocycle in $H^*(U_n)^{U_n}$ of degree $p$. Notice that the real
dimension of $U_n$ is $n^2$. Then $\chi(\omega)\in
HC^{n^2-k}(C^\infty_c(\Gamma \ltimes U_n))$ is defined as
\[
\chi(\omega)(f_0,f_1, \cdots, f_{n^2-p})=\int_{U_n}\omega \wedge
\big(f_0df_1\cdots df_{n^2-k}\big)|_{id},
\]
where $f_0, \cdots, f_{n^2-p}$ are elements in $C_{c}^\infty(\Gamma
\ltimes U_n)$,  $d$ is the de Rham differential on
$\Omega^*(U_n)\rtimes \Gamma$, and $(\cdots )|_{id}$ stands for the
restriction of an element in $\Omega^*(U_n)\rtimes \Gamma$ to the
identity component.

In the following we prove that $\chi(\omega)$ defines an
$(n^2-p)$-trace on  the reduced $C^*$-algebra $C(U_n)\rtimes \Gamma$
and hence the extension exists. Our proof is essentially the same as
Connes' arguments in \cite{connes:transverse}.  However, it
simplifies greatly in that we have an isometric (rather than almost
isometric) action at hand.  In this case the initial step of finding
a Banach sub-algebra $B$ on which to construct an n-trace is
unnecessary, since the map $\lambda: C_{c}^\infty(\Gamma \ltimes
U_n) \to End_{A}(\cale)$ is a
*-homomorphism, hence is bounded.  Thus one can take $B = C(U_n)\rtimes \Gamma$ with the reduced $C^*$-norm.

% Let $V_{1}$ denote the complexified cotangent bundle of $U_{n}$,
% $T^*U_n \otimes \C$, and let
% $V_j$  stand for the exterior algebra bundles $\wedge^j T^*U_n \otimes
% \C$, $j=1, 2, \cdots$. We remark that on each
% $V_j$ there is a hermitian metric $<,>_{j}$  which is
% invariant under the $\Gamma$ action and which induces a metric on the
% space of sections of the bundle $V_j$. Let $r:\Gamma \ltimes U_n\to
% U_n$ be the projection onto the first component, and
% $\cale_{j}=C(U_n, V_j)\rtimes_{red} \Gamma$ be the completion of
% $C_c(\Gamma \ltimes U_n, r^*(V_j))$ respect to the norm $||\cdot
% ||_{j}$ induced by $<,>_{j}$.

% We remark that the left multiplication of $C(U_n)\rtimes \Gamma$ on
% $\cale_{1}$ defines a linear map $\lambda$ from $C_c(U_n\rtimes
% \Gamma)$ to $\End_{A}(\cale_1)$, where $A$ is $C(U_n)\rtimes \Gamma$.
% According to \cite{co:fundamental}[Lemma 3.2], the map $\lambda$ is
% closable. Let $B$ be the domain of the closure of $\lambda$ and equip
% $B$ with the graph norm $|||x|||=\sup(||x||, ||\lambda(x)||)$.  In
% general, $B$ is a dense Banach subalgebra of $A$, but need not be
% closed under holomorphic functional calculus. However, because
% $U^\delta_n$ action on $U_n$ is isometric, by
% \cite{co:fundamental}[Prop. 3.5] $B$ does satisfy this condition and hence
% \[
% K_*(B)\rightarrow K_*(B)
% \]
% is an isomorphism.

What we will show next is that elements in the image, $\hat
\chi(H^*(U_n)^{U^n})$, are $n$-traces which by
\cite{connes:transverse}[Thm. 2.7] will define linear maps from
$K_*(A)$ to $\complex$.

We recall the definition of an $n$-trace on a Banach algebra
$A$. An $n$-trace on $A$ is an $n+1$ linear functional
$\tau$ on a dense subalgebra $\cala$ of $A$ such that
\begin{enumerate}
\item $\tau$ is a cyclic cocycle on $\cala$.
\item for any $a^i\in \cala$, $i=1,\cdots, n$, there exists
$C=C_{a^1, \cdots, a^n}<\infty$ such that:
\begin{equation}\label{eq:estimate}
|\hat{\tau}((x^1da^1)(x^2da^2)\cdots (x^nda^n))|\leq
C||x^1||\cdots ||x^n||, \ \ \ x^i\in A.
\end{equation}
\end{enumerate}

For our purpose, we consider the Banach algebra to be $C(U_n)\rtimes
\Gamma$.  Given any $\omega\in \Omega^*(U_n)^{U_n}$ with
$d\omega=0$, we have that $\chi(\omega)$ defines a cyclic cocycle on
$C_{c}^\infty(\Gamma \ltimes U_n)$.  To prove that $\chi(\omega)$
defines an $m$-trace, it suffices to prove the estimate
(\ref{eq:estimate}) for any $a^i\in \cala$ and $x^i\in A$, $i=1,
\cdots,n$.

We define a convolution product on the following collection of
spaces $\bigoplus_j \cale_j=\bigoplus_j C_c(\Gamma \ltimes U_n,
r^*(\wedge^jT^*_\complex U_n))$  by $*: \cale_i\otimes
\cale_j\rightarrow \cale_{i+j}$,
\[
\phi* \psi(\alpha)=\sum_{\alpha=\beta\gamma}\phi(\beta)\wedge
\beta(\psi(\gamma)), \ \ \ \alpha\in \Gamma \ltimes U_n.
\]

Assume that $\omega\in \Omega^p(U_n)^{U_n}$ with $d\omega=0$. The
following formula $\tr_\omega(\phi)=\int_{U_n}\omega(x)\phi(x,id)$
defines a linear map on $\cale_{n^2-p}$ with the following
property
\[
\tr_\omega(\phi * \psi)=\tr_{\omega}((-1)^{jk}\psi* \phi),\ \ \
\psi\in \cale_j,\ \phi \in \cale_k,\ j+k=n^2-p.
\]

The proof of the above trace property goes as follows.
\[
\begin{split}
&\tr_\omega(\phi * \psi)\\
=&\int_{U_n}\omega(x)\sum_{\alpha\beta=id}\phi(x, \alpha)\wedge
\psi(\alpha(x), \beta)\\
=&\int_{U_n}\omega(x)(-1)^{jk}\sum_{\alpha\beta=id}\psi(\alpha(x),
\beta)\wedge \phi(x, \alpha)\\
\end{split}
\]
\[
\begin{split}
=&\int_{U_n}\omega(\alpha^{-1}(y))(-1)^{jk}\sum_{\alpha\beta=id}\psi(y,\beta)\wedge
\phi(\beta(y), \alpha)\\
=&\int_{U_n}\omega(y)(-1)^{jk}\sum_{\alpha\beta=id}\psi(y,\beta)\wedge
\phi(\beta(y), \alpha)\\
=&(-1)^{jk}\tr_\omega(\psi*\phi).
\end{split}
\]
In the above proof, $y$ is $\alpha(x)$ and we have used that
$\alpha^{-1}=\beta$, and $\omega$ is $\Gamma$ invariant.

We prove the following properties for $\tr_\omega$, which is an
analog of \cite{connes:transverse}[Lemma 4.3].
\begin{lemma}
\noindent
\begin{enumerate}
\item For any $\phi\in \cale_{n^2-p}$, there is a constant
$C_\phi<\infty$ such that for any $f\in C(U_n)\rtimes\Gamma$,
one has
\[
|\tr_\omega(f\phi)|\leq C_\phi||f||_A,
\]
where $||f||_A$ is the $C^*$-algebra norm of the reduced crossed
product $C^*$-algebra $C(U_n)\rtimes \Gamma$.
\item For and $\phi_1, \cdots, \phi_{n^2-p}\in \cale_1$, there
exists a constant $C_{\phi_1, \cdots,\phi_{n^2-p}}<\infty$, such
that for any $f_1, \cdots, f_n\in C(U_n)\rtimes \Gamma$,
\[
\tr_\omega(\phi_1f_1* \phi_2f_2*\cdots* \phi_{n^2-p}f_{n^2-p} )\leq
C_{\phi_1,\cdots, \phi_{n^2-p}}\prod_i ||f_i||.
\]
\end{enumerate}
\end{lemma}
\begin{proof}

(1) We recall that for any $f\in C(U_n)\rtimes \Gamma$, the following estimates of the $C^*$-norm of $||f||_A$ holds, i.e. $\sup_{x\in U_n, \gamma\in \Gamma}|f(x,\gamma)|\leq ||f||_A$.  Now for $\tr_\omega(f\phi)$, we have that
\[
\begin{split}
|\tr_\omega(f\phi)|&=|\int_{U_n}\sum_{\alpha\beta=id}\omega(x)f(x, \alpha)\alpha(\phi(\alpha(x), \beta))|\leq\int_{u_n}|\omega(x)|\sum_{\alpha\beta=id}|f(x,\alpha)||\alpha(\phi(\alpha(x),\beta))|\\
&\leq ||f||_A\int_{U_n}|\omega(x)|\sum_{\alpha}|\alpha(\phi(\alpha(x), \alpha^{-1}))|\leq C_\phi ||f||_A .
\end{split}
\]
We remark that  $\sum_\alpha |\alpha(\phi(\alpha(x), \alpha^{-1}))|$ is a finite sum, and therefore the integral is finite. \\

(2) We apply Connes \cite{connes:transverse}[Thm. 3.7] to the following situation. \\

\noindent{\bf Theorem}: {\em Let $\Psi$ be an $m$-linear function on $\cale_1$ satisfying the following conditions

\indent{(a)} $\Psi(\xi_1, \cdots, \xi_j f, \xi_{j+1}, \cdots, \xi_m)=\Psi(\xi_1, \cdots, f\xi_{j+1}, \cdots, \xi_m)$ for $j=1,\cdots, m-1$, $\xi_k\in \cale_1 $ and $f\in C(U_n)\rtimes \Gamma$.

\indent{(b)} For any $\xi_1, \cdots, \xi_m\in \cale_1$, there exists $C<\infty$ such that
\[
|\Psi(\xi_1, \cdots, \xi_m f)|\leq C||f||_A, \ \ \ \ f\in C(U_n)\rtimes \Gamma.
\]
Then for any $f_0\in C(U_n)\rtimes \Gamma$, and $\xi_1,
\cdots, \xi_m\in \cale_1$, there exists $C'<\infty$ with
\[
|\Psi(f_0f_1\xi_1, f_2\xi_2, \cdots, f_m\xi_m)|\leq
C'||f_1||\cdots||f_m||,\ f_i\in C(U_n)\rtimes \Gamma.
\]
}

To apply the above theorem to our situation. We introduce $\Psi_\omega$ an $n^2-p$ linear functional $\cale_1$ as follows, for $\xi_1, \cdots, \xi_{n^2-p}\in \cale_1$,
\[
\Psi_\omega(\xi_1, \cdots, \xi_{n^2-p})=\int_{U_n}\omega(x)\wedge
\xi_1*\cdots * \xi_{n^1-p}.
\]

It is easy to see that $\Psi_\omega$ satisfies condition (a) in the
above Theorem as $\xi_if* \xi_{i+1}=\xi_i* f* \xi_{i+1}=\xi_i*
f\xi_{i+1}$ for any $i=1, \cdots, n^1-p-1$. For condition (b), we
consider that $\Psi(\xi_1, \cdots,
\xi_{n^1-p}f)=\tr_\omega(\xi_1*\cdots* \xi_{n^2-p}* f )=\tr_\omega(f
\xi_1* \cdots* \xi_{n^2-p})$. As all $\xi_i$ belongs to $\cale_1$,
$\xi_1* \cdots * \xi_{n^2-p}$ again belongs to $\cale_1$. Hence, by
Lemma 3.2 (1), $|\Psi_\omega(\xi_1,\cdots,
\xi_{n^2-p}f)|=|\tr_\omega(f\xi_1*\cdots *\xi_{n^2-p})|\leq
C_{\xi_1, \cdots, \xi_{n^2-p}}||f||_A$. Therefore both condition (a)
and (b) are satisfied, and the above Theorem implies that
\[
|\Psi_\omega(f_0f_1\xi_1, f_2\xi_2, \cdots, f_{n^2-p}\xi_{n^2-p})|\leq C'||f_1||\cdots ||f_{n^2-p}||.
\]

Now using the tracial property, we have that
\[
\begin{split}
&|\tr_\omega(\phi_1f_1*\cdots \phi_{n^2-p}f_{n^2-p})|=|\tr_{\omega}(f_{n^2-p}\phi_1* f_1\phi_2*\cdots* f_{n^2-p-1}\phi_{n^2-p})\\
=&|\Psi_\omega(f_{n^2-p}\phi_1, \cdots, f_{n^2-p-1}\phi_{n^2-p})|\leq C'||f_1||\cdots ||f_{n^2-p}||.
\end{split}
\]
Hereby Lemma 4.2 is proved.
\end{proof}

Now by Lemma 4.2 and the tracial property, we can easily obtain the
condition (\ref{eq:estimate}), and we conclude that for any
$\omega\in \Omega^p(U_n)^{U_n}$ with $d\omega=0$, $\hat
\chi(\omega)\in HC^{n^2-p}(C_{c}^\infty(\Gamma \ltimes U_n))$
defines a linear function from $K_*(C(U_n)\rtimes\Gamma)$ to
$\complex$.

\end{proof}

\begin{remark} We can extend $U_n$ in Theorem \ref{thm:extend} to a
  general manifold $V$ with a vector bundle $E$ which is equipped with
  a almost isometric action of $\Gamma$. Then the similar statements
  hold for the closed $\Gamma$ invariant differential forms on $V$. In
  \cite{connes:transverse}, Connes applies this idea to fundamental
  cocycles. In the context, one can see that he definitely has this
  kind of generalization in mind.
\end{remark}

\begin{remark}We remark that Jiang \cite{jiang} proved that if a
discrete group $\Gamma$ is rapid decay and acts isometrically on a
closed oriented riemannian manifold $V$ preserving the orientation,
then cocycles on $E\Gamma\times _\Gamma V$ of polynomial growth (See
\cite{jiang} for the precise statement) can be paired with the
$K$-theory group of the reduced $C^*$-algebra $C_0(V)\rtimes
\Gamma$. Our results in this section show that any
$\Gamma$-invariant cocycles on $V$ can be paired with the $K$-theory
of the $C^*$-algebra for any orientation preserving isometrical
action of a discrete group $\Gamma$ on an oriented riemannian
manifold. We do not need to assume the group to be rapid decay.
\end{remark}

\section{Transgressed classes and cyclic classes}
The goal of the present section is to study the transgressed
Chern character of a flat trivialized $U_n$ bundle.   By providing a
simplicial construction of the universal transgressed Chern character
we will be able to show that it is in the image of the map $\chi$.
Indeed, the proof actually shows that any class in $\chi
(H^{*}(U_{n})^{U_{n}})$ is obtained by   transgression
 from an invariant polynomial applied to curvature forms.

We start from the following  fibration
\[
\overline{BU}_n\rightarrow BU_n^\delta\rightarrow BU_n,
\]
where $\overline{BU}_n$ is the homotopy fiber of the map from
$BU_n^\delta \rightarrow BU_n$.  We can take $\overline{BU}_n$ to be
the realization of a simplicial space, $EU_n^\delta\times
_{U_n^\delta} U_n$. Now, one knows that $BU_n$ is the classifying
space for $U_n$ bundles, $BU^\delta_n$ is the classifying space for
flat $U_n$ bundles, and $\overline{BU}_n$ is the classifying space
for trivialized flat $U_n$ bundles. In particular, we will see that
on $EU^\delta_n\times_{U_n^\delta}U_n$ there is a  flat $U_n$ bundle
equipped with a canonical trivialization.

Consider
the universal flat $U_n$ bundle over $BU_n^\delta$, which can be
identified with $EU_n^\delta$. Therefore, the associated flat
principal $U_n$ bundle over $BU_n^\delta$ can be identified with
$EU^\delta_n\times _{U_n ^\delta} U_n\rightarrow BU^\delta_n$.
Considering the map of classifying space, we have the following
diagram
\begin{equation}
\label{diagram1}
\begin{diagram}
\node{EU_n^\delta\times_{U_n^\delta}(U_n\times U_n)}\arrow{s,l,T}{\pi'} \arrow{e,t,T}{\iota'} \node{EU_n^\delta\times_{U_n^\delta}U_n}\arrow{s,r,T}{\pi}\\
\node{EU_n^\delta\times_{U_n^\delta}U_n}\arrow{e,t,T}{\iota}\node{BU^\delta_n}
\end{diagram}
\end{equation}
where $EU_n^\delta\times_{U_n^\delta}(U_n\times U_n)$ with $U_n^\delta$ acting on $U_n\times U_n$ diagonally is the pullback respect to the maps $\iota$ and $\pi$.

We write the maps
$\pi, \pi', \iota, \iota'$ in coordinates.  We write a point in $EU_n^\delta\times_{U_n^\delta} U_n$ by $(z,x)$ with $z\in EU_n^\delta$ and $x\in U_n$. Then $\pi(z,x)=\iota(z,x)= [z]$, where $[z]$ stands for the $U_n^\delta$ orbit of $z$ in $EU_n^\delta$. And $(z,x,y)$ with $z\in EU_n^\delta, x, y\in U_n$ form a local $U_n^\delta$-equivariant coordinates for $EU_n^\delta \times_{U_n^\delta}(U_n\times U_n)$ such that $\iota'(z,x,y)=(z,y)$ and $\pi'(z,x,y)=(z,x)$.

On the bundle $\pi':EU_n^\delta\times_{U_n^\delta}(U_n\times U_n)\to EU_n^\delta\times_{U_n^\delta}U_n$, there are two natural connections. One is the pull back of the flat connection on $\pi: EU_n^\delta\times _{U_n^\delta}U_n\to BU_n^\delta$, the other is from a canonical trivialization $EU_n^\delta\times_{U_n^\delta}(U_n\times U_n)\to (EU_n^\delta\times _{U_n^\delta}U_n) \times U_n$.

\begin{enumerate}
\item The flat connection on $\pi: EU_n^\delta\times _{U_n^\delta}U_n\to BU_n^\delta$ is given by Dupont \cite{du:simplicial} in a simplicial way. We consider the universal $U_n^\delta$ bundle $\pi: \overline{NU}_n^\delta\to NU_n^\delta$. $\overline{NU}_n^\delta$ is the simplicial manifold with $\overline{NU}_n^\delta(k)=U_n^\delta\times \cdots \times U_n^\delta$ ($k+1$ copies), where the face operators are defined by leaving out one of the components. Similarly, $NU_n^\delta$ is simplicial manifold with $NU_n^\delta(k)=U_n^\delta\times \cdots \times U_n^\delta$ ($k$-copies) with the face operator equal to the multiplication of consecutive two components.  $EU_n^\delta$ and $BU_n^\delta$ are the corresponding geometrical realizations of $\overline{NU}_n^\delta$ and $NU_n^\delta$.
There is a canonical map from $\overline{NU}_n^\delta\to NU_n^\delta$ defined by
\[
\pi(x_0, \cdots, x_k)=(x_0x_1^{-1}, x_1x_2^{-1}, \cdots, x_{k-1}x_k^{-1}).
\]

We consider the following maps $pr_i:\overline{NU}_n^\delta\times_{U_n^\delta}U_n\to U_n$ by
\[
pr_i(x_0, \cdots, x_k; x)=x_ix.
\]
Furthermore, let $\theta$ be the Maurer-Cartan form on $U_n$, and $(t^0, \cdots, t^k)$ be the barycentric coordinates of the simplex $\Delta_k$. Then $\Theta=\sum_{i}t^i pr^*_i(\theta)$ defines a connection one form on $EU_n^\delta\times_{U_n^\delta}U_n\to BU_n^\delta$. The pullback of this connection through $\iota'$ defines a connection one form on $EU_n^\delta\times_{U^\delta_n}(U_n\times U_n)$, which is denoted by $\Theta_0=\iota'^*(\Theta)$.
\item We construct a canonical trivialization of
$EU_n^\delta\times_{U_n^\delta}(U_n\times U_n)\to
EU_n^\delta\times _{U_n^\delta}U_n$ by using the natural section
$\sigma:EU^\delta_n\times_{U_n^\delta}U_n\to EU_n^\delta\times
_{U_n^\delta}(U_n\times U_n)$ given by $\sigma(z,x)=(z,x,x)$. Accordingly we
have a map $\Sigma:EU_n^\delta\times_{U_n^\delta}(U_n\times
U_n)\to U_n$ defined by $\Sigma(z,x,y)=y^{-1}x$.  If $\theta$ is the
Maurer-Cartan form on $U_n$, then $\Sigma^*(\theta)$ defines a
connection one form on $EU_n^\delta\times_{U_n^\delta}(U_n\times
U_n)$ which is denoted by $\Theta_1=\Sigma^*(\theta)$ and which is
associated to the section $\sigma$.
\end{enumerate}

We consider the following connection $D=s\Theta_0+(1-s)\Theta_1$ on
the space $(EU_n^\delta\times_{U_n^\delta}(U_n\times U_n))\times
[0,1]$, where $s$ is the coordinate on the interval $[0,1]$. Let
$R_s$ be the curvature of $D$.  Let $\mathfrak{U}_n$ be the Lie
algebra of $U_n$. For any $\calp \in
sym(\mathfrak{U}_n^*)^{\mathfrak {U}_n}$, we define $T\calp \in
\Omega^*(EU_n^\delta\times_{U_n^\delta}(U_n\times U_n))$ by
\[
T\calp=\int_0^1 \calp(R_s)
\]

Using the expression $D=s\Theta_0+(1-s)\Theta_1$, we have that
\[
R_s=ds\wedge \Theta_0+sd\Theta_0-ds\wedge
\Theta_1+(1-s)d\Theta_1+\frac{1}{2}s^2\Theta_0^2+(1-s)s\Theta_0\wedge
\Theta_1+ \frac{1}{2}(1-s)^2\Theta_1^2.
\]
It is easy to check that $R_s$ is $\mathfrak{U}_n$ horizontal, and
therefore descends to a $2$-form on
$(EU_n^\delta\times_{U_n^\delta}U_n)\times [0,1]$. Consider the
projection map $\pi':EU_n^\delta\times_{U_n^\delta}(U_n\times
U_n)\times [0,1]\to EU_n^\delta\times_{U_n^\delta}U_n\times
[0,1]$, which  forgets the second  $U_n$ factor. We observe that the restrictions of
$\Theta_0=\sum_i t^ipr_i^*(\theta)$ and $d\Theta_0=\sum_i dt^i\wedge pr_i^*(\theta)+t^ipr_i^*(d\theta)$ to $BU_n^\delta $ vanish as $U_n^\delta$ is equipped with discrete topology and any differential form with a positive degree on a discrete set always vanishes. We conclude
that restrictions of those terms in $R_s$ containing form $\Theta_0$ or
$d\Theta_0$ vanish on the base $EU_n^\delta\times_{U_n^\delta}U_n\times
[0,1]$. Thus  $R_s$ simplifies on
$EU_n^\delta\times_{U_n^\delta}U_n\times [0,1]$ to
\[
-ds\wedge \Theta_1+(1-s)d\Theta_1+\frac{1}{2}(1-s)^2\Theta_1^2.
\]
For example,  when we apply this to the polynomial $\calp (X)=\Tr(\exp(-\frac{1}{2\pi}X))\in
sym(\mathfrak{U}_n^*)^{\mathfrak{U}_n}$,  we obtain a formula for the
transgressed Chern character,
\[
Tch=\int_0^1 \Tr(\exp(R_s))=\int_0^1 ds \sum_k
\frac{1}{k!}\frac{(-1)^{k+1}}{(2\pi)^{k+1}}\Tr(\Theta_1(sd\Theta_1+\frac{1}{2}s^2\Theta_1^2)^k).
\]

Consider next $\Theta_1$ and $d\Theta_1$, which are simplicial
differential forms on the components of the simplicial manifold
$EU_n^\delta\times_{U_n^\delta}(U_n\times U_n)(k)=(\Delta_k\times
\overline{NU}_n^\delta(k))\times_{U^\delta_n} (U_n\times U_n)$
independent of the coordinates on $\Delta_k$. Therefore $Tch$ is a
differential form on the simplicial manifold
$EU_n^\delta\times_{U_n^\delta}U_n$ independent of the $\Delta_k$
components in $(\Delta_k\times
\overline{NU}_n^\delta)\times_{U^\delta_n} U_n$.

Let $I$  be the map from the space of differential forms
$\Omega^l((\Delta_k\times \overline{NU}_n^\delta)\times_{U^\delta_n}
U_n)$ on the simplicial manifold $EU_n^\delta\times_{U_n^\delta}U_n$
to the group cochain complex $C^{k}(U_n^\delta; \Omega^{l-k}(U_n))$
of $\Omega^*(U_n)$ valued $U_n^\delta$ cochains. $I$ is realized by
integration along the $\Delta_k$ component of $(\Delta_k\times
\overline{NU}_n^\delta)\times_{U^\delta_n}  U_n$. Finally, observing
that $T\calp$ is independent of the $\Delta_k$ coordinates, we
conclude that $I(T\calp)$ gives rise to a zero dimensional
$U_n^\delta$ cocycle with values in $\Omega^*(U_n)^{U_n}$.

Finally,  we will relate characteristic numbers obtained from these
classes to secondary classes for trivialized flat bundles.  Note that we can view these constructions in the context of the
following diagram.
\begin{equation}
  \begin{diagram}
    \node{U_{n}} \arrow{s,l}{} \arrow{e,t}{=} \node{U_{n}}
    \arrow{s,r}{}\\
\node{V = \tilde M \x_{\Gamma} U_{n}}\arrow{s,l}{}\arrow{e,t}{\hat \alpha}
\node{EU_{n}^{\delta} \x_{U_{n}^{\delta}} U_{n} =
  \overline{BU}_{n}} \arrow{s,r}{k}\\
\node{M} \arrow{r,t}{\alpha} \arrow{ne,t}{\hat \theta}\node{BU_{n}^{\delta}}\\
\end{diagram}
\end{equation}

 The bundle that was considered in (\ref{diagram1}) is the pull-back of
the left side to its total space.  As a trivialized flat bundle it
is classified by the map $\hat \alpha$. The universal class $T\calp$
is in $H^{*}(\overline{BU}_{n})$.  If we assume in addition the
existence of a trivialization, $\theta :\tilde M \x_{\Gamma} U_{n}
\to M \x U_{n}$, then there is a lift of $\alpha$ determined by
$\theta$ which we will denote by $\hat \theta : M \to
\overline{BU}_{n}$.  The pull back of the universal class by $\hat
\theta$ we will denote by
\begin{equation}
  \hat \theta^{*}(T\calp) = T\calp (\alpha,\theta).
\end{equation}
This will play a role in the index theorems of the next section.

 Combining this with Theorem \ref{thm:extend} we obtain the following result.
\begin{proposition}
\label{tp} The following holds.
\begin{itemize}
\item[1)]  Let $[c] \in \chi(H^{*}(U_{n})^{U_{n}})$.  Then there
exist $i\geq 1$ and
  polynomials $\calp_1, \cdots, \calp_i \in sym(\mathfrak{U}_n^*)^{\mathfrak{U}_n}$
  such that $[c]= \hat \alpha^{*}([T\calp_1]\cdots [T\calp_i])$.
\item[2)] Let $\calp \in
  sym(\mathfrak{U}_n^*)^{\mathfrak{U}_n}$.  Then the associated transgressed
  class, $\hat \alpha^{*}(T\calp)$ is equal to $\chi([c])$ for some cyclic class $[c]$.
\item[3)] The cyclic class from (2)  extends to define a linear map
  on the
  K-theory of the reduced crossed-product $C^{*}$-algebra.
  \begin{equation}
    [c] :K_{*}(C(U_n)\rtimes\Gamma) \to \C.
  \end{equation}
\end{itemize}
\end{proposition}

% Moreover, as a corollary to the above analysis and Theorem
% \ref{thm:extend}, we obtain the following theorem.

% \begin{theorem}\label{thm:trans-chern} For any $\calp(X) \in
%   sym(\mathfrak{U}_n^*)^{\mathfrak{U}_n}$, the associated transgressed
%   class, $T\calp$ is equal to $\chi([c])$ for some cyclic class $[c]$
%   which extends to define a linear map on  K-theory
%   \begin{equation}
%     [c] :K_{*}(C(U_n)\rtimes_{red}\Gamma) \to \C
%   \end{equation}
%   where $\Gamma$ is a subgroup of $U_n^\delta$.
% \end{theorem}

% We relate our universal transgressed Chern character to the one
% defined in \cite{do-he-ka:eta}. What really counts in our
% construction is the following $U_n\times U_n\rightarrow U_n$
% bundle. (We are taking semidirect product of this bundle with
% $EU_n^\delta$ over $U_n^\delta$ action. However, this does not
% contribute to the final answer as we have seen that in the image
% of $I$ only the $U_n$ components show up in the final answer.) On
% this bundle there are two splittings. One is from the natural map
% from $U_n\rightarrow U_n\times U_n$ by embedding $U_n$ into the
% first component, the other is the diagonal embedding
% $U_n\rightarrow U_n\times U_n$. The first connection corresponds
% to $\Theta_0$ in our construction, and the second connection
% corresponds to $\Theta_1$. And then after taking the image under
% $I$ inside $\Omega^*(U_n)^{U_n}$, we see that the
% transgressed Chern character constructed in this section agrees
% with the one constructed in \cite{do-he-ka:eta}[Section 4].

\section{Higher index theorems and secondary classes}

In this section we will apply Connes' index theorem, \cite{connes:book},
% , as adapted to
% smooth etale groupoids in \cite{Lott-Gorokhovsky, Gorokhovsky}
to
obtain explicit formulas for pairing our cyclic classes with index
classes in K-theory of
certain operators.  This will lead in the next section to what appears
to be a new family of classes in the cohomology of a group which will
satisfy the Novikov conjecture.  They are the analog of the Gelfand-Fuchs
classes which Connes considered in \cite{connes:transverse}.

Recall the geometric setting as described in Section 2.
 Let $M$ be a closed $Spin^{c}$
manifold and let $\alpha: \Gamma =\pi_{1}(M) \to U_{n}$ be an
injective homomorphism.  Consider the associated flat, foliated,
principal $U_{n}$ bundle, $V = \tilde M \x_{\alpha} U_{n} \to M$.
Let $\DD$ denote the Dirac operator on $M$ and $\tilde \DD$ its lift
to $\tilde M$.  Then $\tilde \DD$ descends to a leafwise elliptic
operator, $\DD_{\alpha}$, on $V$ with respect to the foliation
$\calf$ with leaves the images of $\tilde M \x \{g\}$ in $V$.

In this context one can apply the higher index theorem of Connes,
\cite{connes:book}.  It provides a topological formula for the
pairing of the index of a leafwise elliptic operator, as an element
of the K-theory, with a cyclic cocycle in
$HP^{*}(C^{\infty}_{c}(U_{n} \rtimes \Gamma))$.  We will use the
form presented by Connes that involves a localization map, $\lambda:
HP^{*}(C^{\infty}_{c}(U_{n} \rtimes \Gamma)) \to H^{*}(V)$,
(\cite{connes:book}, p. 274), which is a left inverse of the map
$\phi: H^{*}(V) \to HP^{*}(C^{\infty}_{c}(U_{n} \rtimes \Gamma))$,
as in (\ref{phi}). We will make strong use of the fact that $\phi$,
and hence also $\lambda$, are isomorphisms.  Suppose also that we
are given a vector bundle $E$ on $V$.

There is an invariant transverse measure for the foliation of $V$
and an associated Ruelle-Sullivan current which we will denote by
$\Lambda$.
\begin{theorem}[Connes]  Let $Ch(\sigma(\DD_{\alpha}))$ denote the Chern
  character of the symbol of $\DD_{\alpha}$, and $Td(T\calf\otimes  \C)$
  the Todd class of the complexified tangent bundle along the leaves.  Let $\Phi$ be the
  Thom isomorphism for $T\calf$ and let $Ind(\DD_{\alpha} \otimes E) \in
  K_{*}( C(U_{n}) \rtimes \Gamma)$ be the K-theoretic index of
  the operator $\DD_{\alpha} \otimes E$.  Let $[c] \in
  HP^{*}(C_{c}^{\infty}(U_{n} \rtimes \Gamma))$ be given. Then one has
  \begin{equation}
\label{index}
    <Ind(\DD_{\alpha} \otimes E), [c]> = <
    \lambda(c)\cup Ch(E) \cup \Phi_{T^{*}\calf}^{-1}(Ch(\sigma(\DD_{\alpha}))\cup
    Td(T\calf\otimes  \C), [\Lambda]>.
\end{equation}
\end{theorem}

We next want to obtain a more precise version of the topological
pairing in the case that the cyclic class $[c] = \hat \chi([u])$,
with $[u] \in H^{*}(U_{n})^{U_{n}}$.  According to
Proposition~\ref{tp} there are invariant polynomials $\calp_1,
\cdots, \calp_i$ so that $\chi([u]) = \hat \alpha^{*}
(T\calp_1\cdots T\calp_i)$. Thus, (\ref{index}) can be rewritten as
\begin{equation}
  <Ind(\DD_{\alpha} \otimes E), [c]> = <
    \hat \alpha^{*} (T\calp_1\cdots T\calp_i) \cup p^{*}Ch(E) \cup \Phi_{T^{*}\calf}^{-1}(Ch(\sigma(\DD_{\alpha}))\cup
    Td(T\calf\otimes  \C), [\Lambda]>.
\end{equation}
If we further assume the existence of a trivialization $\theta: \tilde
M \x_{\Gamma} U_{n} \to M \x U_{n}$, then we can simplify to
\begin{equation}
  <Ind(\DD_{\alpha} \otimes E), [c]> = <
     \hat \pi^{*} \hat \theta^{*} (T\calp_1\cdots T\calp_i) \cup p^{*}Ch(E) \cup \Phi_{T^{*}\calf}^{-1}(Ch(\sigma(\DD_{\alpha}))\cup
    Td(T\calf\otimes  \C), [\Lambda]>,
\end{equation}

Recall that we defined $T\calp(\theta,\alpha) = \hat
\theta^{*}(T\calp_1\cdots T\calp_i)$. Following the argument in
\cite{do-he-ka:eta}, this may be used to obtain the index by a
pairing on the base $M$. Thus, we finally obtain the desired
formula.
\begin{proposition}
\label{charnumber} Let $V = \tilde M \x_{\alpha} U_{n}$ be the
foliated, flat, bundle obtained from the representation
$\alpha:\pi_{1}(M) \to U_{n}$.  Let $\theta : \tilde M \x_{\alpha}
U_{n} \to M \x U_{n}$ be a trivialization.  Let $[c] \in
\chi(H^{*}(U_{n})^{U_{n}})$ correspond to $T\calp_1\cdots T\calp_i$.
Then we have
\begin{equation}
  <Ind(\DD_{\alpha} \otimes E), [c]> = <
       T\calp(\theta,\alpha) \cup Ch(E) \cup \Phi_{T^{*}M}^{-1}(Ch(\sigma(\DD))\cup
    Td(TM\otimes  \C), [M]>.
\end{equation}
\end{proposition}
If we now specialize to the case where the bundle $E$ is the spinors,
so that $\DD \otimes E$ is the signature operator on $M$, which we
denote $D^{sign}_{M}$. Then the operator on $\tilde M \x_{\Gamma}
U_{n}$ is the leafwise signature operator, and we obtain
the following formula which we will apply in the next section.

\begin{proposition}
\label{indextheorem}
\begin{equation}
  <Ind(D^{sign}_{M,\alpha}), [c]> = <
       T\calp(\theta,\alpha) \cup L(M), [M]>.
\end{equation}
\end{proposition}

\section{Homotopy invariance of characteristic numbers}

In this section we will obtain homotopy invariance results for
characteristic numbers obtained from the cyclic classes considered
above.  We will describe two methods for doing this.  The first is
based on Connes' method of extending cocycles,
c.f. \cite{connes:book, connes:transverse}.    The second was
communicated to us by Guoliang Yu and we will present his argument as
an alternative approach.

\subsection{First method}
\label{first}
Let $h:M_{1} \to M_{2}$ be an orientation preserving homotopy
equivalence of closed, oriented manifolds.  Let $\alpha_{i}:
\Gamma_{i} = \pi_{1}(M_{i}) \to U_{n}$ be unitary representations such that
$\alpha_{2}h_{*} = \alpha_{1}$.

The first step is to do a suspension
operation.  For this we lift $h$  to a homotopy equivalence
between universal covers, $\tilde h : \tilde M_{1} \to \tilde M_{2}$,
and this descends to a leafwise homotopy equivalence between the
foliated, flat, principal bundles,
\begin{equation}
  \hat h : V_{1} = \tilde M_{1} \x_{\alpha_{1}} U_{n} \to \tilde M_{2}
  \x_{\alpha_{2}} U_{n} = V_{2}.
\end{equation}

The signature operators on $M_{i}$  induce operators on $V_{i}$,
$D^{Sign}_{M_{i},\alpha_{i}}$ which are the leafwise signature
operators. The map $\hat h$ yields a Morita equivalence which
induces a map on the K-theory of the associated (\'etale) foliation
algebras,
\begin{equation}
  \hat h_{*}: K_{*}( C(U_{n}) \rtimes_{\alpha_{1}} \Gamma_{1})) \to K_{*}(C(U_{n}) \rtimes_{\alpha_{2}} \Gamma_{2})),
\end{equation}
and it follows from \cite{K-M,hilsum-skandalis} that
\begin{equation}
\label{homotopyinvariance}
\hat h_{*}(\Ind
(D^{Sign}_{M_{1},\alpha_{1}})) = \Ind (D^{Sign}_{M_{2},\alpha_{2}}).
\end{equation}
Recall that one can always represent the K-theory class of the index
in the image of $K_{*}(C_{c}^{\infty}(\Gamma\ltimes U_{n}))$ in such
a way that any of our cocycles can pair with it.  However, the
equivalence in (\ref{homotopyinvariance}) is in $K_{*}(C(U_{n})
\rtimes \Gamma)$, so there is no guarantee that $<\hat h_{*}(\Ind
(D^{Sign}_{M_{1},\alpha_{1}})), [c]>$ will agree with $<\Ind
  (D^{Sign}_{M_{1},\alpha_{1}}), \hat h^{*}([c])>$.  If the cocycle
  $[c]$ extends to the crossed-product $C^{*}$-algebra, then we do
  obtain equality.
Thus, we have the following version of the homotopy invariance
property expressed in the Novikov conjecture.

\begin{proposition}
\begin{equation}
\label{equality}
  <\Ind (D^{Sign}_{M_{2},\alpha_{2}}), [c]> =  <\hat h_{*}(\Ind (D^{Sign}_{M_{1},\alpha_{1}})), [c]> = <\Ind
  (D^{Sign}_{M_{1},\alpha_{1}}), \hat h^{*}([c])>,
\end{equation}
where $\hat h^{*}$ is the map induced on cyclic theory of the
crossed-products.
\end{proposition}
Now, if we interpret the pairings using the results from the
previous sections, these equalities can be expressed in terms of
characteristic numbers.  Suppose that the cyclic class $[c]$
corresponds to a cohomology class in $\chi(H^{*}(U_{n})^{U_{n}})$.
There is an integer $N$ so that the bundle defined by
$\alpha_{1}^{N}$ is trivial, and we will choose an explicit
trivialization $\theta$.  With this data, we can apply the index
formulas (\ref{indextheorem}) to (\ref{equality}) and obtain
homotopy invariance results for certain characteristic numbers.

% we can, after bringing in a
% trivialization, express it as a transgressed class,
% $T\calp(\alpha_{1}^{N},\theta)$.

\begin{theorem}
  $<h^{*}(L(M_{2}))\cup T\calp (\alpha_{1}^{N},\theta), [M_{1}]> = <L(M_{2})\cup
  T\calp (\alpha_{2}^{N}, \theta), [M_{2}]>.$
\end{theorem}
\begin{proof}
  By the index theorem, (\ref{index}), the first and third terms in
  (\ref{equality}) can be rewritten in terms of topological pairings
  with the Ruelle-Sullivan currents associated to the Haar measure,
  and then, because of the form of $\chi([c])$, the pairing can be
  done on the base.  The resulting formula is as stated.
\end{proof}

This should be viewed as a homotopy invariance property for the local
part of an invariant which is  defined in special cases.  For the
case of $Tch$, the transgressed Chern character, the invariant is the
$\rho$-invariant of Atiyah-Patodi-Singer.  Homotopy invariance for the
$\rho$-invariant has been proved under the condition that the Baum-Connes conjecture
holds for  $C^{*}_{max}(\Gamma)$.  However, it has been shown to be true in $\R / \Q$
in general, \cite{farber-levine, keswani, piazza-schick}.

\subsection{Second method}
\label{second}
The second approach, communicated to us by Guoliang Yu, does not require
the  extending of cocycles.  We will sketch the argument, which is an
adaptation  of the usual proof that injectivity of the assembly map
implies the Novikov conjecture.

\begin{proposition}
  The assembly map
  \begin{equation}
    \mu : K^{top}(U_{n}^{\delta};C(U_{n})) \to K_{*}(C(U_{n})\rtimes_{max}
    U_{n}^{\delta})
  \end{equation}
is injective.  Moreover, there is an injective map
\begin{equation}
  \xi : K_{*}(\underline{EU^{\delta}_{n}}\x_{U_{n}^{\delta}} U_{n})
  \otimes \Q = K_{*}(\overline{BU}_{n}) \otimes \Q
  \to K^{top}(U_{n}^{\delta};C(U_{n}))\otimes \Q
\end{equation}

\end{proposition}

Assuming the same data as in Section \ref{first}, $(M_{i},
\alpha_{i},\theta_{i})$,  we get leafwise signature operators and a
map $\hat \theta_{i}: M_{i} \to \overline{BU}_{n}$ with the property
that $\mu(\xi[(M_{i},\theta_{i})]) =
Ind(D^{Sign}_{M_{i},\alpha_{i}})$.  According to
\cite{hilsum-skandalis}, one has $Ind(D^{Sign}_{M_{2},\alpha_{2}}) =
\hat h_{*}(Ind(D^{Sign}_{M_{1},\alpha_{1}}))$ and since $\mu$ and
$\xi$ are
injective we have that $[(M_{2},\hat\theta_{2}h)] =
[(M_{1},\hat \theta_{1} )]$.  Our goal is to show that for any class
 $[u] \in HC^{*}(\calh(U_{n}^{\delta}\ltimes U_n))_{n^2}$  we have
 \begin{equation}
   <Ind(D^{Sign}_{M_{2},\alpha_{2}}), \chi([u])> =
<\hat h_{*}(Ind(D^{Sign}_{M_{1},\alpha_{1}})), \chi([u])>.
 \end{equation}
Since the map $\phi : \oplus_k H^{*+2k}(\overline{BU}_{n}) \to
HP^{*-n^2}(C_{c}^{\infty}( U_{n}^{\delta}\ltimes U_n))$ is an
isomorphism, there is a class $x \in H^{*}(\overline{BU}_{n})$
satisfying
\begin{equation}
  <Ch([M_{i},\theta_{i}]),x> = <Ind(D^{Sign}_{M_{i},\alpha_{i}}), \chi([u])>.
\end{equation}

Since the left side is independent of $i$ so is the right, and we obtain,
\begin{proposition}

\begin{equation}
   <Ind(D^{Sign}_{M_{2},\alpha_{2}}), \chi([u])> =
<(Ind(D^{Sign}_{M_{1},\alpha_{1}})), \hat h^{*}\chi([u])>.
 \end{equation}
\end{proposition}

Note that this last step is the reverse of the usual procedure.  Here
we have a cyclic class and we must find a cohomology class
corresponding to it, while in Connes' argument one has a cohomology
class and one must find a corresponding cyclic class which also extends to
K-theory of the crossed-product.

\section{Concluding remarks}

For a finitely presented group, $\Gamma$, let $LP(\Gamma) \subseteq
H^{*}(\Gamma,\Q)$ denote the elements obtained by the above process.
That is, given a finite dimensional unitary representation of
$\Gamma$, we consider $\chi(HC^{*}(\calh(U_{n} \rtimes
U_{n}^{\delta})_{n^2}) \subseteq H^{*-n^2}(\tilde M \x_{\Gamma}
U_{n})$ and we integrate these classes along the fiber to obtain a
subgroup of $H^{*}(B\Gamma)$.  By the above arguments, pairing these
classes with the index of  signature operators gives homotopy
invariant characteristic numbers.  This is the unitary analog of
Connes' work on Gelfand-Fuchs classes, (LP stands for
Lazarov-Pasternack classes).

\bibliographystyle{alpha}

\end{document}